\documentclass[a4paper,11pt]{article}
\usepackage{amsfonts}
\usepackage{amsmath}%
\usepackage{amsthm}
\usepackage{amsfonts}%
\usepackage{amssymb}%
\usepackage{graphicx}
\usepackage{graphics}
\usepackage[latin1]{inputenc}
\usepackage{hyperref}
\newtheorem{prop}{Proposition}
\newtheorem{theorem}{Theorem}
\newtheorem{rem}{Remark}
\newtheorem{lem}{Lemma}
\newtheorem{defi}{Definition}

\newtheorem{cor}{Corollary}

\title{Invariant bipartite random graphs on $\mathbb{R}^d$}
\author{ Fabio Marcellus Lopes }

\begin{document}

\maketitle

\begin{abstract}
\small{Suppose that red and blue points occur in $\mathbb{R}^d$ according to two simple point process  with finite intensities $\lambda_{\mathcal{R}}$ and $\lambda_{\mathcal{B}}$, respectively. Furthermore, let $\nu$ and $\mu$ be two probability distributions on the strictly positive integers. Assign independently a random number of stubs (half-edges) to each red and blue point with laws $\nu$ and $\mu$, respectively. We are interested in translation-invariant schemes to match stubs between points of different colors in order to obtain random bipartite graphs in which each point has a  prescribed degree distribution with law $\nu$ or $\mu$ depending on its color. Let $X$ and $Y$ be random variables with law $\nu$ and $\mu$, respectively. For a large class of point processes we show that we can obtain such translation-invariant schemes matching a.s. all stubs if and only if
\[   \lambda_{\mathcal{R}} \mathbb{E}(X)= \lambda_{\mathcal{B}} \mathbb{E}(Y), \]
allowing $\infty$ in both sides, when both laws have infinite mean. Furthermore, we study a particular scheme based on the Gale-Shapley stable marriage \cite{gale}. For this scheme we give sufficient conditions on $X$ and $Y$ for the presence and absence of infinite components. These results are  two-color versions of those obtained in \cite{Mia2}.}
\end{abstract}

% Suppose that we have red and blue points in $\mathbb{R}^{d}$, where the set of points of each color was obtained from different point process in  $\mathbb{R}^{d}$  with finite positive intensity. Let each point be independently equipped with a random number of stubs (half-edges) with a commom positive integer-valued law for the points of the same color. We are interested in translation-invariant schemes to match stubs between points of different colors in order to obtain a graph in which each component is a simple bipartite graph and the vertex degree has the prescribed law for each color.\\

% We show that it is possible if and only if  For a particular matching scheme based on the Gale-Shapley stable marriage we show for a large class of point processes that it yields almost surely a unique graph with the desired properties for a given configuration of points and stubs, if both laws for the number of stubs have finite mean and a simple condition relating the mean values and intensities is satisfied. Furthermore, in the case where the points are from two independent homogeneous Poisson processes, we give sufficient conditions on the laws for the presence or absence of an infinite component. Our work extends to the two color case the main results obtained in \cite{Mia2}, where this model with a unique Poisson process in $\mathbb{R}^d$ was considered.

\section{Introduction}

Let $\mathcal{R}$ and $\mathcal{B}$ (\textbf{red and blue points})  be two translation-invariant simple point processes on $\mathbb{R}^d$, jointly ergodic under translations, with finite intensities $\lambda_{\mathcal{R}}$ and $\lambda_{\mathcal{B}}$, respectively. Furthermore, let $\nu$ and $\mu$ be two probability laws on the strictly positive integers. We assign independently to each red and blue point a random number of stubs (half-edges) with law $\nu$ and $\mu$, respectively. Our aim is to study schemes for pairing the stubs in order to obtain translation-invariant simple bipartite random graphs whose vertices are points of $\mathcal{R}$ and $\mathcal{B}$, where the degree of each vertex has law $\nu$ or $\mu$ depending on its color, and where edges between pairs of points of the same color are not allowed. The first natural question is how different the two point processes and their stub laws can be for schemes matching a.s. all stubs to be possible. This question was first asked in \cite{Mia}. In the simplest case when we have a.s. one stub per point, it is easy to believe that the intensities of the point processes must coincide in order to obtain a perfect matching between them. Indeed, in \cite{Yuval1} this claim is proved and many other properties of so-called \textbf{two-color perfect matchings} between two point processes are studied, in particular, bounds on the matchings distances are given. In our model, with multiple stubs per point, new questions arise.  For example, can we always give translation-invariant schemes to pair the stubs that a.s. yield at least one infinite component? Can we give schemes that a.s. give only finite components? For pairing schemes which can lead to both kinds of components, can we give conditions on $\nu$ and $\mu$ that guarantee percolation and non-percolation, respectively? In \cite{Mia2} these questions were studied  for the one-color case, that is, one Poisson process with i.i.d assigned stubs to each point, and the pairing schemes, of course, allow connections between points of the same color. Particular attention is spent on a matching scheme based on the Gale-Shapley stable marriage \cite{gale}. In this work we give a sharp condition for when a matching of a.s. all stubs in the two-color case exists and we extend some of the main results of \cite{Mia2} to the two-color case.\\
\indent  Next, we describe more formally our problem and the random objects we will work with.
 The \textbf{support} (or point-set) of $\mathcal{R}$ is the random set $[\mathcal{R}]:= \{ x \in \mathbb{R}^d : \mathcal{R}(\{x\})>0\}$, its points are called \textbf{red points}. Analogously, we write $[\mathcal{B}]$ for the point-set of the process $\mathcal{B}$, and its points are called \textbf{blue points}. In general, for any random point measure $\Lambda$ we write $[\Lambda]$ for its support. The \textbf{intensity} of a translation-invariant point process is the expected number of points in a Euclidean ball of unit volume.\\
\indent  Let  $X$ and $Y$ be random variables with law $\nu$ and $\mu$, respectively, and let $\eta_{\mathcal{R}}$ be a random integer-valued measure on $\mathbb{R}^d$ with the same support as $\mathcal{R}$, which conditionally on $\mathcal{R}$, assigns i.i.d values with law $X$ to the elements of $[\mathcal{R}]$. Similarly, for $\mathcal{B}$, let $\eta_{\mathcal{B}}$ assign i.i.d values with law $Y$ to the elements of $[\mathcal{B}]$.  The pairs $(\mathcal{R},\eta_{\mathcal{R}} )$ and $(\mathcal{B} ,\eta_{\mathcal{B}} )$ are two marked point processes with positive integer-valued marks. For $x \in [\mathcal{R}]$, we write $X_{x}$ for $\eta_{\mathcal{R}}(\{x\})$ and, for $y \in [\mathcal{B}]$, we write $Y_{y}$ for $\eta_{\mathcal{B}}(\{y\})$, which we interpret as the number of stubs at the red point $x$ and the number of stubs at the blue point $y$, respectively. Sometimes we refer to the stubs as red or blue depending on the color of the point to which they are assigned. For a marked point process, we refer to the expected number of stubs in a Euclidean ball of unit volume as the \textbf{stub intensity}.\\
\indent A \textbf{two-color multi-matching scheme} for two marked processes $(\mathcal{R},\eta_{\mathcal{R}} )$ and $(\mathcal{B} ,\eta_{\mathcal{B}} )$ is a point process on the space of unordered pairs of points in $\mathbb{R}^d$ with the pro\-per\-ty that almost surely for every unordered pair $(x,y)\in [\mathcal{M}]$ we have $x\in [\mathcal{R}]$ and $y\in [\mathcal{B}]$, and such that in the bipartite graph $B=B(\mathcal{R},\mathcal{B},\mathcal{M})$ with vertex set $[\mathcal{R}]\cup[\mathcal{B}]$ and edge set $[\mathcal{M}]$, each vertex $x \in [\mathcal{R}]$ and each vertex $y \in [\mathcal{B}]$ has degree $X_{x}$ and $Y_{y}$, respectively. The two-color perfect matching mentioned before is the special case in which all points of both processes have degree one almost surely. If a vertex $x \in [\mathcal{R}]$ (or $y \in [\mathcal{B}]$) has degree at most $X_{x}$ (or $Y_{y}$), we talk about a \textbf{two-color partial multi-matching scheme}. We only consider simple  two-color multi-matchings and partial multi-matchings schemes where the bipartite graph B has no self-loops and no multiple edges, and that are translation-invariant, meaning that $\mathcal{M}$ is invariant in law under the action of all translations of $\mathbb{R}^d$. Let $\mathbb{P}$ be the probability measure governing $\left( \mathcal{B}, \eta_{\mathcal{R}},\mathcal{B} ,\eta_{\mathcal{B}}, \mathcal{M} \right)$. We say that a two-color partial multi-matching is a \textbf{factor} if $\mathcal{M}$ is a deterministic function of $\left( \mathcal{B}, \eta_{\mathcal{R}},\mathcal{B} ,\eta_{\mathcal{B}}, \mathcal{M} \right)$, that is, if it does not involve any extra randomness. We introduce the \textbf{Palm process} $\left( \mathcal{B}^{*}, \eta_{\mathcal{R}}^{*},\mathcal{B}^{*} ,\eta_{\mathcal{B}}^{*}, \mathcal{M}^{*} \right)$, with law $\mathbb{P}^{*}$ and expectation $\mathbb{E}^{*}$, in which we condition on the presence of a red point at the origin, while taking the mark processes, the pairing scheme and $\mathcal{B}$ as a stationary background. See e.g \cite{kalle, daley2} for details about Palm processes and point processes. For the Palm version of our process, we denote by $C$ the volume of the component of the red point at the origin, that is, the number of vertices that can be reached by a path in $B(\mathcal{R}^{*}, \mathcal{B}^{*}, \mathcal{M}^{*})$.\\
\indent Similarly, a \textbf{one-color multi-matching} of $(\mathcal{R},\eta_{\mathcal{R}})$ is a point processes $\mathcal{M}$ on the space of unordered pairs of point in $\mathbb{R}^d$, with the property that a.s. for every pair $(x,y)\in \mathcal{M}$ we have that $x,y\in \mathcal{R}$, and such that in the graph $G=G(\mathcal{R},\mathcal{M})$ with vertex set $[\mathcal{R}]$ and edge set $[\mathcal{M}]$ , each vertex $x$ has degree $X_{x}$.\\
\indent  Let $|.|$ denote  Euclidean distance. A set $S \subset \mathbb{R}^d$ is said to be \textbf{non-equidistant} if there are no distinct points $x,y,u,w \in S$ with $|x-y|=|x-z|$ or $|x-y|=|u-w|$. A \textbf{descending chain} is an infinite sequence $\{ x_{i}\}_{i\geq1} \subset S$ such that $|x_{i} - x_{i-1}|$ is strictly decreasing. Most of our results require that the underlying point processes have support that a.s. is non-equidistant and has no descending chains.  We observe that both conditions hold for Poisson processes, for a proof see \cite{hagg}. Also, if $\mathcal{B}$ and $\mathcal{R}$ are two independent Poisson processes, then almost surely $[\mathcal{R}]\cup [\mathcal{B}]$ is non-equidistant and has no descending chains.\\
\indent Our first result gives a necessary and sufficient condition for the existence of a two-color multi-matching scheme for two marked point processes $(\mathcal{R},\eta_{\mathcal{R}})$ and $(\mathcal{B},\eta_{\mathcal{B}})$. 
\begin{theorem}
\label{teo1}
\indent Let $(\mathcal{R},\eta_{\mathcal{R}})$ and $(\mathcal{B},\eta_{\mathcal{B}})$ be two marked point processes on $\mathbb{R}^d$, jointly ergodic under translations, with simple ground processes of finite intensities $\lambda_{\mathcal{R}}$ and $\lambda_{\mathcal{B}}$, and with i.i.d positive integer-valued marks with laws $X$ and $Y$, respectively. Suppose that almost surely $[\mathcal{R}]\cup[\mathcal{B}]$ is non-equidistant, and has no descending chains. 
Then there exists a two-color multi-matching scheme for $(\mathcal{R},\eta_{\mathcal{R}})$ and $(\mathcal{B},\eta_{\mathcal{B}})$ if and only if
\begin{equation}
\label{condition}
 \lambda_{\mathcal{R}} \mathbb{E}(X)= \lambda_{\mathcal{B}} \mathbb{E}(Y).
\end{equation}
In (\ref{condition}) we also allow $\infty$ on both sides, meaning that both processes have marks with infinite mean.
\end{theorem}
\indent The only-if part follows from a simple application of the so-called mass transport principle (see Lemma \ref{mass}). For the if part we introduce two procedures leading to two-color multi-matchings, one that applies when both sides in (\ref{condition}) are finite and  another one for the infinite case. These two procedures are related to the notion of a stable matching. 
\begin{defi}
A two-color multi-matching scheme $\mathcal{M}$ is said to be a \textbf{two-color stable multi-matching} if a.s., for any pair of points $x\in [\mathcal{R}]$ and $y\in [\mathcal{B}]$, either they are linked by an edge or at least of them has no incident edges longer than $|x-y|$.
\end{defi} 
Next, we state the main results on the percolation questions from the beginning of the introduction. As already mentioned they extend results proved by \cite{Mia2} for the one-color case. Our proofs are based on similar constructions combined with ideas of \cite{Yuval1} on two-color matchings. Now, we leave the setting of more general point processes and work with two independent marked Poisson processes.
\begin{theorem}
\label{schemes}
Let $(\mathcal{R},\eta_{\mathcal{R}})$ and $(\mathcal{B},\eta_{\mathcal{B}})$ be two independent marked Poisson processes on $\mathbb{R}^d$, $d\geq1$, with finite intensities $\lambda_{\mathcal{R}}$ and $\lambda_{\mathcal{B}}$, respectively, and satisfying the condition $(\ref{condition})$. 
\begin{itemize}
\item [(a)]If  $(\mathcal{R},\eta_{\mathcal{R}})$ and $(\mathcal{B},\eta_{\mathcal{B}})$ have the same law, we refer to such situation as  \textbf{ the symmetric case}. Then there exists a simple translation-invariant factor matching such that $\mathbb{P}^{*}(C<\infty)=1$.
\item [(b)]  If $\mathbb{P}(Y\geq2)\geq \mathbb{P}(X\geq2)>0$, then there exists a simple translation-invariant factor matching scheme with $\mathbb{P}^{*}(C=\infty \mid X_{0}\geq2)=1$.
\end{itemize}
\end{theorem}

The next result gives sufficient conditions on the degree distributions $X$ and $Y$ that guarantee the existence and non-existence, respectively of a component with infinitely many vertices for the two-color stable multi-matching of two independent marked Poisson processes. 
\begin{theorem}
\label{perc}
 Let $(\mathcal{R},\eta_{\mathcal{R}} )$ and $(\mathcal{B} ,\eta_{\mathcal{B}} )$ be two marked Poisson processes on $\mathbb{R}^d$, jointly ergodic under translations, with finite intensities $\lambda_{\mathcal{R}}$ and $\lambda_{\mathcal{B}}$, and i.i.d marks with laws $X$ and $Y$ both with finite mean, and satisfying (\ref{condition}).
 Consider the two color stable multi-matching.
\begin{itemize}
 \item [(a)]  For any $d \geq 2$, there exists a $k=k(d)$ such that if $\mathbb{P}\left( Y\geq k  \right)=\mathbb{P}\left( X\geq k  \right)=1$, then $\mathbb{P}^{*} \left( C=\infty \right)>0$.
\item [(b)] For any $d\geq1$, we have that if $\mathbb{P}\left( Y\leq 2  \right)=\mathbb{P}\left( X\leq 2  \right)=1$, and $\mathbb{P}\left( Y= 1  \right)>0$, then $\mathbb{P}^{*} \left( C=\infty \right)=0.$
\end{itemize}
\end{theorem}
\indent It can be seen that the two sufficient conditions are quite apart from each other. It is an open problem to determine sharp conditions for our model as well as for the one color case. Recently, for the one-color case some advances have been made in $d=1$ with $\mathbb{P}(X=2)=1$, see \cite{Mia5}, specifically, rigorous arguments are provided that show that there exists an infinite component if a certain event in a finite interval has large enough probability. The latter is supported by simulations.\\
\indent Invariant spatial random graphs have been studied recently, for example, \cite{Mia4, Mia3, jona} study automorphism invariant random graphs on lattices with prescribed degree distribution, \cite{Ferrari, Yuval2, Yuval1} study translation-invariant trees and matchings on point processes in $\mathbb{R}^d$, and  \cite{Mia, Mia2} study translation-invariant random graphs for Poisson processes in $\mathbb{R}^d$ with prescribed degree distribution.\\
\indent In the next section we introduce an algorithm that yields a two-color stable multi-matching and prove uniqueness for the resulting matching. Furthermore, we prove Theorem \ref{teo1}. In the third section, we prove Theorems \ref{schemes} and \ref{perc}. Finally, in the last section, we pose some open questions.
\section{Two-color stable multi-matching } 
\indent Next, we introduce an iterative procedure that generates a two-color stable multi-matching for the two marked point processes $(\mathcal{R},\eta_{\mathcal{R}} )$ and $(\mathcal{B} ,\eta_{\mathcal{B}} )$. Our prodedure is a two-color extension of the algorithm proposed in \cite{Mia} for one-color multi-matchings. That algorithm, in turn, is a multi-matching generalization of the iterated mutually closest matching algorithm from \cite{Yuval1}. We will refer to our algorithm as \textbf{2CIMC}.\\
\indent  Following \cite{Yuval1}, we call a pair of points $x$, $y$ \textbf{potential partners} if one is red while the other is blue. We say that two potential partners $x$ and $y$ are \textbf{mutually closest} if $y$ is the closest potential partner to $x$ and $x$ is the closest potential partner to $y$. We denote by $\mathcal{M}(x)$ a partner of a point $x$; that is  $\{x,\mathcal{M}(x)\}\in [\mathcal{M}]$. Next we describe the 2CIMC algorithm.\\ 
\indent Given the point configurations $[\mathcal{R}]$ and $[\mathcal{B}]$, we start by creating an edge between each mutually closest potential partners in $[\mathcal{R}]\cup[\mathcal{B}]$, and then removing one stub from each of these points. In the next step, we consider the set of points which still have at least one stub after the previous step. We call two potential partners compatible if no edge was created between them in the previous step. An edge is created between each compatible mutually closest  potential partners in this set, and again we remove one stub from each of these points. Then the algorithm is iterated.\\ 
\indent By construction, the algorithm above yields almost surely a two-color partial multi-matching that avoids multiple edges and self-loops. The next result shows that if the two processes have the same finite stub intensity then almost surely all stubs are matched and we obtain a two-color stable multi-matching. In the case where the stub intensities are different we obtain almost surely a two-color partial stable multi-matching that exhausts all stubs in the process with smaller stub intensity.
% Step 1. Consider the set $[\mathcal{R}]\cup[\mathcal{B}]$ of all blue and red points. An edge is created between each mutually closest potential partners in the set, and one stub is removed from each of these  points.\\
% 
% Step 2. Consider the set of all points that still have at least one stub after step 1. Two potential partners are called compatible if no edge was created between them in step 1. An edge is created between each compatible mutually closest potential partners in this set, and one stub is removed from each of these points.\\
% \begin{center}
% \vdots \end{center}
% Step \textit{n}. Consider the set of all points that still have at least one stub. Two potential partners are called compatible if no edge has beeb created between them. An edge is created between each compatible mutually closest potential partners in this set, and one stub is removed from each of these points.\\
% \begin{center}
% \vdots
% \end{center}
\begin{prop}
\label{st1}
Let $(\mathcal{R},\eta_{\mathcal{R}} )$ and $(\mathcal{B} ,\eta_{\mathcal{B}} )$ be two marked Point processes on $\mathbb{R}^d$, jointly ergodic under translations, with ground processes with finite intensities $\lambda_{\mathcal{R}}$ and $\lambda_{\mathcal{B}}$, and i.i.d marks with law $X$ and $Y$ both with finite mean, respectively.
Suppose that almost surely $[\mathcal{R}]\cup[\mathcal{B}]$ is non-equidistant and has no descending chains.\\
\begin{itemize}
\item [(a)] If (\ref{condition}) holds, then almost surely the 2CIMC algorithm described above exhausts the set of stubs, and the limiting graph (after an infinite number of iterations) is a two color stable multi-matching. No other two color stable multi-matching of $(\mathcal{R},\eta_{\mathcal{R}} )$ and $(\mathcal{B} ,\eta_{\mathcal{B}} )$ exists.
\item [(b)]If (\ref{condition}) holds with a strict inequality, the 2CIMC algorithm yields a translation-invariant two-color partial  stable multi-matching scheme that almost surely exhausts all stubs in the process with smaller stub intensity. 
\end{itemize}
\end{prop}
\begin{rem}
\label{stable}
For the case when $\nu(\{1\})=\mu(\{1\})=1$, the result is an application of the two-color stable matching of \cite[Proposition 9]{Yuval1}.
\end{rem}
\begin{rem}
\label{infinite1}
We note that the 2CIMC algorithm can be applied when some pairs of vertices already have an edge between them and additional connections between such vertices are prohibited. When the existing edges form a translation-invariant process the same argument as in the proof of  Proposition \ref{st1} shows that the procedure yields a  two-color partial multi-matching that a.s. exhausts all stubs of the process with lower stub intensity.
\end{rem}
We state without proof the following well-known lemma that is usually called the mass transport principle.
\begin{lem}(\textbf{Mass transport principle}, \cite[Lemma 8 (ii)]{Yuval1}).
\label{mass}
 Suppose $T$ is a random non-negative measure on $\mathbb{R}^d \times \mathbb{R}^d $ such that $T(A,B):=T(A\times B)$ and $T(A+w,B+w)$ are equal in law for all $w \in \mathbb{Z}^d $. Then
\[  \mathbb{E}T(Q,\mathbb{R}^d ) =\mathbb{E}T(\mathbb{R}^d ,Q). \] 
\end{lem}
\indent Given a translation-invariant two-color partial multi-matching, the point process of matched red stubs is the point process that counts the number of unordered pairs $\{x,\mathcal{M}(x) \}\subset \mathbb{R}^d$ where $x \in [\mathcal{R}]$. By definition, the intensity of this process is the expected number of red stubs matched by $\mathcal{M}$ in a fixed set of unit volume, we also call this intensity the \textbf{matched red stub intensity}. The process of the remaining red stubs (associated with the given two-color partial multi-matching) is a marked point process that assigns to each point  $x\in [\mathcal{R}]$ an integer which is the difference between its initial number of stubs $X_{x}$ and the number of partners $x$ has in the two-color partial multi-matching. If we define the \textbf{unmatched red stub intensity} to be the expected number of unmatched red stubs in a fixed set of unit volume, we have that this number is the difference between the initial stub intensity of the red points and the matched red stub intensity. These two point processes are defined on the same probability space supporting  $(\mathcal{R},\eta_{\mathcal{R}},\mathcal{B} ,\eta_{\mathcal{B}},\mathcal{M})$. We define similar objects for the blue points.
\begin{lem}{(Fairness)}
\label{fair}
Let $(\mathcal{R},\eta_{\mathcal{R}} )$ and $(\mathcal{B} ,\eta_{\mathcal{B}} )$ be two marked point processes on $\mathbb{R}^d$ with ground processes with finite intensities $\lambda_{\mathcal{R}}$ and $\lambda_{\mathcal{B}}$, and i.i.d marks with law $\eta_{\mathcal{R}}$ and $\eta_{\mathcal{B}}$, respectively. 
Let $\mathcal{M}$ be a translation-invariant two-color partial multi-matching of the two marked processes. Then the point processes of matched red stubs and of matched blue stubs have equal matched stub intensity.
\end{lem}
\indent \textit{Proof of Lemma \ref{fair}}: Apply Lemma \ref{mass} to the mass transport
% \[ T(A,B):= \#\{ (x,\mathcal{M}(x))_{\mbox{\tiny{unordered}}}\in [\mathcal{M}] : x \in [\mathcal{R}]\cap A ; \mathcal{M}(x) \in [\mathcal{B}]\cap B \}\]
in which each red point in A with matched stubs sends unit mass to each one of its partners in B. Then
$\mathbb{E}T(Q, \mathbb{R}^d )$  is the intensity of matched red stubs, while $\mathbb{E}T(\mathbb{R}^d , Q)$ is the intensity of matched blue stubs. $\square$\\
\indent This lemma is a multi-matching analogue of  Lemma 7 in \cite{Yuval1}, and it gives an immediate corollary.
\begin{cor}
\label{corr}
If there exists a translation-invariant two-color multi-matching for $(\mathcal{R},\eta_{\mathcal{R}} )$ and $(\mathcal{B} ,\eta_{\mathcal{B}} )$, then $\lambda_{\mathcal{R}}\mathbb{E}(X)$ and $\lambda_{\mathcal{B}}\mathbb{E}(Y)$ must be equal.
\end{cor}
\indent \textit{Proof of Proposition \ref{st1} (a)}: Let $\mathcal{N}_{\mathcal{R}}$ and $\mathcal{N}_{\mathcal{B}}$ be the point processes of the red and blue points, respectively, with at least one unmatched stub after the 2CIMC algorithm is completed. Then  $\mathcal{N}_{\mathcal{R}}$ and $\mathcal{N}_{\mathcal{B}}$ are ergodic point processes, and hence we have three possible cases:
\begin{itemize}
 \item [(i)] both have a.s infinitely many points,
\item [(ii)] one has a.s infinitely many points and the other has a.s no points,
\item [(iii)] both have a.s no points.
\end{itemize}
Our aim is to rule out the first two cases.\\
\indent (i) To rule out the first case we adapt the argument of the proof of \cite[Proposition 2.2]{Mia2} to the two-color case. First, we call a unordered pair of points  $x \in \mathcal{N}_{\mathcal{R}}$ and   $ y \in \mathcal{N}_{\mathcal{R}}$ compatible if they do not have an edge between them in the configuration obtained after applying the 2CIMC procedure. Then we create a bipartite directed graph $G$ with vertex sets $[\mathcal{N}_{\mathcal{R}} ]$ and $[\mathcal{N}_{\mathcal{B}} ]$ by drawing a directed edge from each point in $[\mathcal{N}_{\mathcal{R}} ]$ to its closest compatible point in $[\mathcal{N}_{\mathcal{B}} ]$, and vice-versa. The closest compatible point exists for each point since the initial number of stubs per point is a.s. finite.\\
\indent The finite components of the graph G would form directed cycles of even sizes (given the bipartite nature), but these cycles cannot be of size greater than two as noted in \cite{Mia2}. So, the finite components must be precisely the cycles of length two.  This however corresponds to two mutually closest compatible points with no edges between them and an unmatched stub at each point, which is impossible since an edge would have been created between them at some stage of the 2CIMC procedure. Hence $G$ has no finite components. This implies that if both $[\mathcal{N}_{\mathcal{R}} ]$ and $[\mathcal{N}_{\mathcal{B}} ]$ are non-empty, there are only infinite components. However, these can only assume the form of infinite descending chains with alternating red and blue points of $[\mathcal{N}_{\mathcal{R}} ]$ and $[\mathcal{N}_{\mathcal{B}} ]$, and, by assumption $[\mathcal{R}]\cup [\mathcal{B}]$ has a.s. no descending chains.\\  
\indent (ii) To rule out the second case, that is, that $[\mathcal{N}_{\mathcal{R}} ]$ has a.s infinitely many points and $[\mathcal{N}_{\mathcal{B}} ]$ is a.s empty, and vice-versa, we adapt the argument of the proof of \cite[Proposition 9]{Yuval1} that shows a similar result for two-color matchings. The key is to use the Fairness Lemma \ref{fair} which says that the processes of matched red stubs and matched blue stubs have equal stub intensity.\\
\indent If we assume that both processes $(\mathcal{R},\eta_{\mathcal{R}} )$ and $(\mathcal{B} ,\eta_{\mathcal{B}})$ satisfy (\ref{condition}), then Lemma \ref{fair} implies that after the 2CIMC algorithm the point processes of remaining red stubs and blue stubs have equal unmatched stub intensity. Since it is not the case if one of them has a.s infinitely many points with unmatched stubs and the other is a.s empty, this case is ruled out.\\
\indent We can conclude that both $[\mathcal{N}_{\mathcal{R}} ]$ and $[\mathcal{N}_{\mathcal{B}} ]$ have a.s no points, which means that the 2CIMC  procedure exhausts a.s. all stubs.\\
\indent That the resulting two-color multi-matching is stable follows from the definition: any unstable pair of compatible points would have had an edge created between them at some step of the two-color multi-matching procedure. The uniqueness of the two-color stable multi-matchings follows by induction over the steps in the algorithm to show that each edge that is present in the resulting configuration must be present in any two-color stable multi-matching for this given configuration of points and stubs.\\ 
\indent \textit{Proof of Proposition \ref{st1} (b)}:
Assume without loss of generality that the red points have a strictly larger stub intensity than the blue points. We want to show that after the procedure we will have $\mathcal{N}_{\mathcal{B}}$ empty almost surely, in other words, all points with remaining stubs are red. To show that we need to rule out the following cases:
\begin{itemize}
 \item [(i')] both $[\mathcal{N}_{\mathcal{R}} ]$ and $[\mathcal{N}_{\mathcal{B}} ]$  have a.s infinitely many points,
\item [(ii')] $\mathcal{N}_{\mathcal{B}}$ has infinitely many points and $\mathcal{N}_{\mathcal{R}}$ has no points a.s.
\end{itemize}
The item (i') is ruled out by the same argument as in the item (i) of part (a) of the proposition. To rule out item (ii') we apply the Fairness Lemma. It implies that the processes of matched red and matched blue stubs have the same matched stub intensity. Since, in (ii'), $\mathcal{N}_{\mathcal{R}}$ is empty almost surely, it would imply that the matched blue stub intensity is  $\lambda_{\mathcal{R}} \mathbb{E}(X)$, but this is strictly bigger than the stub intensity of the blue points. $\square$\\
\indent Our proof of Proposition \ref{st1}(a) does not include the case where both marked point processes have mark processes with infinite mean, since we have used the finiteness of $(\ref{condition})$ in the item (ii) to obtain that the stub intensity of points with remaining stubs is the same in both processes. Below we describe a matching scheme that yields a two-color multi-matching when both point processes have mark processes with infinite mean. Briefly, this matching scheme uses the 2CIMC procedure repeatedly. In each stage the number of stubs that are alllowed to be matched per point is truncated in such a way that at least one of the colors have all their allowed stubs matched. Then we proceed by alternating the color that has all allowed stubs matched at each stage.
\begin{prop}
\label{st2}
Let $(\mathcal{R},\eta_{\mathcal{R}} )$ and $(\mathcal{B} ,\eta_{\mathcal{B}} )$ be two marked point processes on $\mathbb{R}^d$, jointly ergodic under translations, with ground processes with finite intensities $\lambda_{\mathcal{R}}$ and $\lambda_{\mathcal{B}}$, and i.i.d marks with law $X$ and $Y$ both with infinite mean. Suppose that  $[\mathcal{R}]\cup[\mathcal{B}]$ is almost surely non-equidistant and has no descending chains. Then there exists a factor two-color matching scheme that exhausts the set of red and blue stubs almost surely.
\end{prop}
\indent \textit{Proof of Proposition \ref{st2}}:
The idea is to apply the 2CIMC infinitely many times. Each time, we truncate the number of stubs that are allowed to be used for each color in such a way that at least all points of one of the colors have all their allowed stubs matched. And, for the next round we make sure that we alternate the color which has all allowed stubs matched. \\
\indent Since both $X$ and $Y$ have infinite mean and take values on the positive integers, it is possible to choose two increasing sequences of truncated r.v.'s $\{X_{i} \}_{i\geq1} $ and $\{Y_{i}\}_{i\geq1}$ such that $X_{i} \uparrow X$ and $Y_{i} \uparrow Y$, and such that for the odd indices $i$, we have $\lambda_{\mathcal{B}}  \mathbb{E}(Y_{i}) \geq \lambda_{\mathcal{R}}\mathbb{E}(X_{i})$, and for the even indices we have $\lambda_{\mathcal{B}}  \mathbb{E}(Y_{i}) \leq \lambda_{\mathcal{R}}\mathbb{E}(X_{i})$.\\
\indent By the $i-th$ truncated version of the random variable $X$, we mean an integer valued random variable $X_{i}$ that distributes its probability mass equal to $X$ from $1$ to an integer value $J_{i}-1$, and puts the remaining probability mass at the integer value $J_{i}$. By an increasing sequence  of truncated r.v.'s $\{X_{i} \}_{i\geq1}$, we mean that $1\leq J_{1}<J_{2}<\ldots$.  The $i-th$ truncated version of the random variable $Y_{i}$ is defined similarly, and we denote by $K_{i}$ the maximum integer value that it assumes with a positive probability.\\
\indent In the first step, we allow at most $J_{1}$ stubs of each red point to be matched, and similarly at most $K_{1}$ stubs for each blue point, that is, we ignore the other stubs in this first round if they exist. Then, by Proposition \ref{st1}, we obtain a translation-invariant partial two-color multi-matching that exhausts almost surely all allowed stubs on the red points. In other words, in this step we matched all stubs at red points with less than $J_{1}$ stubs and $J_{1}$ stubs at other red points. For the blue points the number of matched stubs per point ranges from $0$ to $K_{1}$, and in addition, by Lemma \ref{fair}, we have that the process of matched red stubs and matched blue stubs have the same intensity.\\
\indent  In the next step, we use the next truncation, and allow at most $J_{2}$ stubs for the red points, including those which were allowed in the previous step, and similarly for the blue points using $K_{2}$ instead. Since in the end of the previous step, we had that the process of matched red stubs and matched blue stubs have the same stub intensity, the inequality $\lambda_{\mathcal{B}}  \mathbb{E}(Y_{2}) \leq \lambda_{\mathcal{R}}\mathbb{E}(X_{2})$ implies that the stub intensity of unmatched allowed blue stubs is smaller than or equal to the stub intensity of unmatched allowed red stubs in this step. We use the 2CIMC procedure with the restriction that pairs of red and blue points that have been connected by an edge in the previous steps cannot be linked again. Since the existing edges were created in a translation-invariant way, by Remark \ref{infinite1}, our procedure also yields a translation-invariant partial two-color multi-matching that exhausts all allowed stubs of the blue points which were not matched in the previous steps. So, we obtain that all stubs at blue points with less than $K_{2}$ stubs and $K_{2}$ stubs at other blue points are matched by the end of the second step.\\
\indent We repeat these steps infinitely many times  using the alternating truncations and the restricted version of the 2CIMC procedure to match pairs of red and blue points with unmatched allowed stubs that did not share an edge created in the previous steps. At the end of each step one of the colors have all stubs matched up to a certain level. Since each point has almost surely a finite number of stubs, for each point at some step all their stubs will be matched. As a result our procedure exhausts almost surely all stubs after infinitely many stages.$\square$\\
\indent \textit{Proof of Theorem \ref{teo1}}: For the only-if part we use  Corollary \ref{corr}. For the if part we have to separate the claim in two cases depending on if (\ref{condition}) is finite or not. In the first case, we obtain the result from Proposition \ref{st1} (a), and for the other one we use Proposition \ref{st2}. $\square$

\section{Percolation for the Poisson case}

In this section, we prove Theorems \ref{schemes} and \ref{perc}.

\subsection{Percolating and non-percolating schemes}

\indent Next, we describe two factor schemes, one that yields a.s. only finite components, and another one which gives a.s. at least one infinite component.\\
\indent \textit{Proof of Theorem \ref{schemes} (a)}: Let $\mathcal{R}_{n}$ denote the process of red points $x\in [\mathcal{R}]$ such that $X_{x}=n$ (i.e, the red points with $n$ stubs), similarly $\mathcal{B}_{n}$ denotes the process of blue points with exactly $n$ stubs. The idea is to partition $[\mathcal{R}_{n}]\cup[\mathcal{B}_{n}]$ into groups of size $2n$ where each group has $n$ points of each color. Then the configuration is taken to consist of bipartite complete graphs on each of these groups.\\
\indent Take $n$ such that $\mathcal{R}_{n}$ and $\mathcal{B}_{n}$ are non-empty, this is possible since $\mu=\nu$. First, to partition $[\mathcal{R}_{n}]$, we assign each red point in $[\mathcal{R}_{n}]$ a type $i \in \{1,2,\ldots,n \}$ as follows. Let $D_{\mathcal{R}_{n}}^{*}$ denote the distance from the origin to the closest other point in the Palm version of $\mathcal{R}_{n}$, and let $0=d_{0},d_{1},\ldots, d_{n-1} ,d_{n}=\infty$ be such that 
\[ \mathbb{P}^{*}(d_{i-1}< D_{\mathcal{R}_{n}}^{*} \leq d_{i}) = \frac{1}{n}, \hspace{0.5cm} i=1,\ldots,n. \]
For $x \in [\mathcal{R}_{n}]$, let $D_{\mathcal{R}_{n}}(x)$ denote the distance to the nearest other point in $[\mathcal{R}_{n}]$. We assign $x \in [\mathcal{R}_{n}]$ type $i$ if $d_{i-1}< D_{\mathcal{R}_{n}}(x) \leq d_{i}$, and let $\mathcal{R}_{n}^{i}$ be the process of points of $\mathcal{R}_{n}$ of type $i$. Since we are in the symmetric case we can use the same numbers $\{d_{0},d_{1},\ldots,d_{n}\}$ to partition $\mathcal{B}_{n}$. Analogously, for $y \in [\mathcal{B}_{n}]$, we define  $D_{\mathcal{B}_{n}}(y)$ and assign the type $i$ if $d_{i-1}< D_{\mathcal{B}_{n}}(y) \leq d_{i}$. Note that the assignment of types does not envolve any extra randomness, and that for each $n$, the processes $\mathcal{R}_{n}^{1},\mathcal{B}_{n}^{1},\ldots,\mathcal{R}_{n}^{n},\mathcal{B}_{n}^{n}$ have equal intensities and are jointly ergodic under translations. 
Since all processes have the same intensity we can use the two-color stable matching repeatedly to construct the groups of size $2n$ as follows. First, for each color, we use the two-color stable matching to match each type $i$ point to a unique type $i+1$ point, for $i=1,2,\ldots,n-1$. The union of these matchings partitions $[\mathcal{R}_{n}]\cup[\mathcal{B}_{n}]$ into monocromatic groups of size $n$. To form groups of size $2n$ with $n$ points of each color we assign to each red point of type $1$ a blue point of type $1$, again using the two color stable matching. 
The components of the graph are bipartite complete graphs formed with each one of these groups.  $\square$
\begin{rem}
The same idea can be applied to cover some asymmetric cases satisfying (\ref{condition}), however, we do not have a general construction for the asymmetric case. 
\end{rem}
% \begin{rem}
% a possible construction for general point processes which are non equidistant and have no descending chains, a factor  or a randomized depending if it is possible to assign types without extra randomness.
% \end{rem}
In order to prove part (b) of Theorem \ref{schemes} we make use of the following result of \cite{Mia2} in which the proof is only sketched since it is part of the proof of Theorem 1 of \cite{Yuval2}. We briefly mention how it is proved because it will be useful to our argument.
\begin{lem}(Lemma 3.1 of \cite{Mia2})\\
For a Poisson process $\mathcal{P}$ with exactly $2$ stubs on each point, there exists a factor matching scheme in which $G$ has a single component consisting of a doubly infinite path.
\end{lem}
\indent The main step to obtain this result is to construct, in a translation-invariant way, a one-ended tree whose vertex set is $[\mathcal{P}]$. Once such one-ended tree has been constructed a doubly infinite path is obtained by first ordering the children of each vertex according to the distance from its parent, and  then ordering all  vertices according to a depth-first search. By linking each pair of vertices that fall next to each other in this ordering by an edge we obtain the desired doubly infinite path (See Theorem 1 in \cite{Yuval2} for details).\\
\indent \textit{Proof of Theorem \ref{schemes} (b)}: First, we will prove the claim for the symmetric case. Let $\mathcal{R}_{\geq2}$ denote the process of red points $x \in [\mathcal{R}]$ with $X_{x}\geq2$, and define $\mathcal{B}_{\geq2}$ similarly. Since both processes have the same intensity, there is a two-color perfect matching of them. We proceed by constructing a translation-invariant one-ended tree whose vertex set is $[\mathcal{B}_{\geq2}]$ as was mentioned before. To obtain a doubly infinite path, however, instead of linking blue vertices that fall next to each other in the ordering, we link each blue vertex to the red vertex matched to the next blue vertex in the ordering by an edge. When this is done we have used two stubs per point and we are left with a doubly infinite path with alternating colors involving all points of $[\mathcal{R}_{\geq2}]\cup[\mathcal{B}_{\geq2}]$. In order to match the points of $[\mathcal{R}]\cup[\mathcal{B}]$ with remaining stubs we apply the 2CIMC procedure with the restriction that we do not allow connections between points that already have an edge between them arising from the connections along the doubly infinite path. Since, by Lemma \ref{fair}, the processes of red and blue points with remaining stubs have equal stub intensity, and the edges created before form a translation-invariant process, we obtain from Remark \ref{infinite1} and Proposition \ref{st1} that the result is a perfect two-color multi-matching.\\
\indent For the asymmetric case, we proceed similarly but first we identify the marked point process with the lowest intensity of points with degree greater than or equal to 2. Then we obtain a two-color partial matching between such points and the points with at least degree two from the other process. Such partial matching must a.s. assign one partner to each point in the process with the lowest intensity of points with degree greater than 2. Then we use such matched pairs to construct the bi-infinite path as we did for the symmetric case. The rest of the proof is equal to the symmetric case. $\square$
\subsection{Percolation for the stable multi-matching}
\indent  Next, we show that the sufficient conditions given in \cite{Mia2} for the existence and absence of an infinite component for the one color stable multi-matching on a Poisson process with i.i.d. degrees can be extend to our model. Our proof of Theorem \ref{perc} is a modification of the proof of \cite[Theorem 1.2]{Mia2}.\\
\indent The proof of Theorem \ref{perc} (a) uses a renormalization argument and a theorem from \cite{schon} concerning domination of r-\textbf{dependent} random fields by product measures.  A random field $\{ X_z \}_{z\in\mathbb{Z}^d}$ is said to be r-dependent if for any two sets $A,B\in \mathbb{Z}^d $ at distance $l_{\infty}$-distance at least $r$ from each other we have that $\{ X_z \}_{z\in A}$ is independent of $\{ X_z \}_{z\in B}$. Below we state without proof the version of the domination result we need.
\begin{theorem}\label{schon1}(\cite{schon}). For each $d\geq2$ and $r\geq1$ there exists a $p_{c}=p_{c}(d,r)<1$ such that the following holds. For any r-dependent random field $\{ X_z \}_{z\in\mathbb{Z}^d}$ satisfying $\mathbb{P}\left(X_{z}=1 \right)=1-\mathbb{P}\left(X_{z}=0 \right)\geq p$ with $p>p_{c}$ the $1$'s in $\{ X_z \}_{z\in\mathbb{Z}^d}$ percolate almost surely. 
\end{theorem}
\indent \textit{Proof of Theorem  \ref{perc} (a)}: For clearness we repeat the renormalization procedure of \cite{Mia2} adapted to our need. We partition $\mathbb{R}^d$ into cubes and declare a cube to be \textbf{good} if it contains at least one point of each color (red and blue) but not too many points of each color and if the same holds for all cubes close to it. The idea is to deduce from Theorem \ref{schon1} that the good cubes can be made to percolate, and we observe that, if each point (blue or red) has sufficiently many stubs attached to it, then each red (blue)  point in a good cube must be connected to each blue (red) point in its adjacent cubes in the two-color stable multi-matching.\\
\indent Let us describe the partition. For $a\in \mathbb{R}$, let $a\mathbb{Z}^d = \{az : z\in \mathbb{Z}^d\}$. We partition $\mathbb{R}^d$ into cubes $\{C_{az}\}_{z\in \mathbb{Z}^d}$ centered at the points of $a\mathbb{Z}^d$ and with side $a$. We call two cubes $C_{az}$ and $C_{ay}$ adjacent if $|z-y|=1$, and write $m=m(d)$ for the smallest integer such that the maximal possible Euclidean distance between points in adjacent cubes does not exceed $ma$. For each cube $C_{az}$ a \textbf{super-cube} $S_{az}$ is defined, consisting of the cube itself along with all cubes $C_{az}$ with $y$ at $l_{\infty}$-distance at most $2m$ from $z$. As a result, a super-cube contains $(4m+1)^d$ cubes.\\
\indent We say that a cube $C_{az}$ is \textbf{acceptable} if it contains at least one point of each color (red and blue) and at most $n=n(d)$ points of each color. The number $n$ will be specified below. A  cube $C_{az}$ is said to be \textbf{good} if all cubes in $S_{az}$ are acceptable. From the independence between $(\mathcal{R},\eta_{\mathcal{R}} )$ and $(\mathcal{B} ,\eta_{\mathcal{B}} )$ we have that the probability that a cube $C_{az}$ is acceptable is equal to
\[
 (1- \mathbb{P}(\mathcal{R}(C_{az}) = 0) - \mathbb{P}(\mathcal{R}(C_{az})>n))(1- \mathbb{P}(\mathcal{B}(C_{az})=0)-  \mathbb{P}(\mathcal{B}(C_{az})>n)). \]
% \begin{equation*}
%  1-\mathbb{P}\left( \mathcal{R}(C_{az})=0 \right)\mathbb{P}\left( \mathcal{B}(C_{az})=0 \right) - \mathbb{P}\left( \mathcal{R}(C_{az})=0 \right)\mathbb{P}\left( \mathcal{B}(C_{az})>0 \right)
% \end{equation*}
% \begin{equation*}
%  - \mathbb{P}\left( \mathcal{R}(C_{az})>0 \right)\mathbb{P}\left( \mathcal{B}(C_{az})=0 \right) - \mathbb{P}\left( \mathcal{R}(C_{az})>n \right)\mathbb{P}\left( \mathcal{B}(C_{az})\leq n \right)  
% \end{equation*}
% \begin{equation*}
%  -\mathbb{P}\left( \mathcal{R}(C_{az})\leq n \right)\mathbb{P}\left( \mathcal{B}(C_{az})> n \right) - \mathbb{P}\left( \mathcal{R}(C_{az})> n \right)\mathbb{P}\left( \mathcal{B}(C_{az}) > n \right)  
% \end{equation*}
% \begin{equation*}
% \geq  1-\mathbb{P}\left( \mathcal{R}(C_{az})=0 \right)\mathbb{P}\left( \mathcal{B}(C_{az})=0 \right) - \mathbb{P}\left( \mathcal{R}(C_{az})=0 \right) - \mathbb{P}\left( \mathcal{B}(C_{az})=0 \right)
% \end{equation*}
% \begin{equation*}
% - \mathbb{P}\left( \mathcal{R}(C_{az})>n \right) - \mathbb{P}\left( \mathcal{B}(C_{az})>n \right) - \mathbb{P}\left( \mathcal{R}(C_{az})>n \right)\mathbb{P}\left( \mathcal{B}(C_{az})>n \right)
% \end{equation*}
\indent By choosing $a$ sufficiently large we can make the probability of having no red (blue) points in $C_{az}$  arbitrarily small, then for a fixed $a$ we can choose $n$ large in order to make the probability of having more than $n$ red (blue) points  in $C_{az}$ arbitrarily small. Hence, the probability that $C_{az}$ is acceptable can be made arbitrarily  close to 1, and consequently, also the probability that $C_{az}$ is good can be made arbitrarily close to 1. In particular, we can make it large enough, as required by Theorem \ref{schon1}, in order to guarantee that the good cubes percolate. Fix such values of $a$ and $n$, and let $k=n(4m+1)^d$ and assume that  $\mathbb{P}\left( Y\geq k  \right)=\mathbb{P}\left( X\geq k  \right)=1$. Now, the same argument of \cite{Mia2}[end of the proof of Theorem 1.2(a)] can be applied for the red and blue points in a good cube to show that they desire all points of different color in the adjacent cubes.\\
% \indent We say that a red (blue) point $x\in [\mathcal{R}]$ ($y\in [\mathcal{B}]$) with $X_{x}$ ($Y_{y}$) stubs \textbf{desires} a blue (red) point $y\in [\mathcal{B}]$ ($x\in [\mathcal{R}]$) if y (x) is one of the $X_{x}$ ($Y_{y}$) nearest points to $x$ in $[\mathcal{B}]$ (to $y$ in $[\mathcal{R}]$).\\
% \indent Then a point $x \in \mathcal{R}$ ($y \in [\mathcal{B}]$) in a good cube $C_{az}$ desires all points in the adjacent cubes: For any point $w \in [\mathcal{B}]$ ($w \in \mathcal{R}$) in an adjacent cube, the distance between $x$ ($y$) and $w$ is at most $ma$, and the Euclidean ball $B_{x}(ma)$ ($B_{y}(ma)$) with radius $ma$ centered at $x$ ($y$) is contained in the supercube $S_{az}$, which contains at most $k$ points of $[\mathcal{B}]$ ($[\mathcal{R}]$).\\
% \indent Since $X_{x}\geq k$ ($Y_{y}\geq k$), it follows that $x$ ($y$) desires all points in $B_{x}(ma)$ ($B_{y}(ma)$), in particular $x$ ($y$) desires $w$. Furthermore, each point $w\in [\mathcal{B}]$ ($w\in [\mathcal{R}]$) in a cube that is adjacent to a good cube $C_{az}$ desires each point in the good cube: Since the distance between $x$ ($y$) and $w$ is at most $ma$, we have that $B_{w}(ma) \subset B_{x}(2ma)$ ($B_{w}(ma) \subset B_{y}(2ma)$). Moreover, the ball $B_{x}(2ma)$ ($B_{y}(2ma)$) is contained in the supercube $S_{az}$, which contains at most $k$ points. Hence $B_{w}(ma)$ contains at most k points, and since $Y_{w}\geq k$ ($X_{w}\geq k$), it follows that $w$ desires all red (blue) points in $B_{w}(ma)$, in particular it desires $x$ ($y$).\\
\indent Since each red point in a good cube desires each blue point in the adjacent cubes and vice-versa, all that remains is to note that two points that desire each other will indeed be matched. This follows from the definition of the two-color stable multi-matching. Hence all points in a good cube are connected to all points of the other color in its own cube and in the adjacent cubes. Since the good cubes percolate this proves the claim. $\square$ \\
\indent Next, we prove Theorem \ref{perc} (b). A modification of the argument in \cite{Mia2}[Theorem 1.2(b)] allows us to prove the claim under the assumption that only one of the processes has a strictly positive probability of degree 1 (if both processes have strictly positive probability of degree 1, it follows from the same argument as in \cite{Mia2}).\\
\indent The next lemma states that the only infinite components that can appear in the graph $B$ obtained from a translation-invariant multi-matching in which all points have degree at most two are a.s. bi-infinite paths.
\begin{lem}
In any translation-invariant two-color multi-matching scheme, a.s $B$ has no component consisting of a singly infinite path.
\label{lemM} 
\end{lem}
\indent \textit{Proof of Lemma \ref{lemM}}: The proof for the two-color case is identical to the one-color case which can be found in \cite[Lemma 5.1]{Mia2}. The proof is obtained by an application of the mass transport principle.$\square$\\
\indent If we restrict the number of stubs per point to be at most $2$, then by the above lemma we obtain that a.s. the only infinite components that can appear in $B$ are bi-infinite paths.

\indent \textit{Proof of Theorem  \ref{perc} }(b): As mentioned above the only infinite components that can appear are a.s. bi-infinite paths. Assume for contradiction that such bi-infinite paths occur with positive probability. For each configuration that contains such a path we will describe a coupled configuration in which, with positive probability, one such path is cut apart into two singly infinite paths. This conflicts with Lemma \ref{lemM}.\\
\indent For any configuration with at least one bi-infinite path, let  $\{x_{i}\}_{i=-\infty}^{\infty}$ be the bi-infinite path with the nearest vertex to the origin, and write $r$ for the second closest red point of the origin on such  bi-infinite path, and $b_{1}$ and $b_{2}$, respectively, for its two blue neighbors on the path. This path will be cut apart by removing the red point $r$ and by re-randomizing the degree of the blue points in the coupled configuration as follows.\\
\indent  Let $\{ Y_{y}  \}_{y \in [\mathcal{B}]}$ be the degree process associated with $(\mathcal{B},\eta_{\mathcal{B}})$. First, construct a coupled configuration in which we introduce for the blue points a modified degree process $\{ \tilde{Y}_{y}  \}_{y \in [\mathcal{B}]}$ with the same law as in the original configuration, and remove the red point $r$. The modified degree process $\{ Y_{y}  \}_{y \in [\mathcal{B}]}$ is obtained as follows. For each $y\in [\mathcal{B}]$, we let $Y_{y}=\tilde{Y}_{y}$ with probability $1-e^{-|y|}$, and with the remaining probability we independently generate a new number of stubs with law $Y$. We call the points which have received a newly generated degree in the modified configuration \textbf{re-randomized points}. We observe that the same Borel-Cantelli argument used in \cite{Mia2} shows that the number of re-randomized point in the modified configuration is almost surely finite:
\[ \mathbb{E}\sum_{z\in[\mathcal{B}] }\mbox{\textbf{1}[z is re-randomized]} = \int_{\mathbb{R}^{d}}e^{-|z|}dz<\infty.\]
\indent Second, we use the fact that the law of $(\mathcal{R}-\delta_{r}, \eta_{\mathcal{R}-\delta_{r}})$ is absolutely continuous $(\prec)$ with respect  to the law of  $(\mathcal{R}, \eta_{\mathcal{R}})$, which is a straightforward modification of \cite{Yuval1}[Lemma 18 (ii)]. By \cite{terry}[Lemma 15], there is a Borel set $S$ with finite Lebesgue measure such that $\mathbb{P}(\mathcal{R}|_{S}=\delta_{r})>0$, where $\mathcal{R}|_{S}$ denotes the restriction of $\mathcal{R}$ to $S$.  By \cite{terry}[Theorem 1 (iii) and Remark 1] and the fact that $(\mathcal{R}-\delta_{r}, \eta_{\mathcal{R}-\delta_{r}})\prec (\mathcal{R}, \eta_{\mathcal{R}})$, we have that $(\mathcal{R}|_{S^c},\eta_{\mathcal{R}|_{S^c}})\prec (\mathcal{R}, \eta_{\mathcal{R}})$; and since the two marked Poisson processes are independent we obtain that $(\mathcal{R}|_{S^c},\eta_{\mathcal{R}|_{S^c}},\mathcal{B},\eta_{\mathcal{B}}) \prec (\mathcal{R}, \eta_{\mathcal{R}},\mathcal{B},\eta_{\mathcal{B}})$.\\
\indent Let $A$ be the event that in the coupled configuration the only re-randomzied points are exactly $b_{1}$ and $b_{2}$, and that both have received a newly generated degree equal to 1. By the finiteness of the set of re-randomized blue points, this event has a positive probability.\\
\indent Now, we claim that, under $A$, the two-color stable multi-matching obtained in the coupled configuration is equal to the one obtained in the original configuration except for the removed red point $r$ and its incident edges $(r,b_{1})$ and $(r,b_{2})$ that do not exist. It is clear that in this matching we have cut apart the bi-infinite path $\{x_{i}\}_{i=-\infty}^{\infty}$  into two singly infinite paths which contradicts Lemma \ref{lemM}.\\
\indent We prove our claim in 2 steps. First, write $B$ for the resulting graph obtained from the original configuration and write $\tilde{B}$ for the graph obtained by removing the point $r$ and its incident edges from $B$. We claim that $\tilde{B}$ is a two-color stable multi-matching for the configurations in $A$. Suppose for contradiction that it is unstable, then there is at least one pair of points $(x,y)$ of different colors in $\tilde{B}$ such that both points have edges which are longer than $|x-y|$. We note that for the points which are still in $\tilde{B}$ there are only two possibilities. Either they are linked to $r$ in $B$ and lost their edges to $r$ in  $\tilde{B}$, or they still have all their partners and edges as in $B$. Since we have not created any new edge or added new pairs of points, the pair of edges which are longer than $|x-y|$ that makes $\tilde{B}$ unstable must also be present in $B$. Furthermore, the pair $(x,y)$ must be stable in $B$. Some thought reveals that the only way it could have been stable in $B$ but not in $\tilde{B}$ is if one of the points had an edge in $B$ which was erased. But that would mean that one of them is $r$ which is not in $\tilde{B}$.
 Second, in our coupled configuration, under $A$, all points have all stubs matched if we match them according to $\tilde{B}$. Since $\tilde{B}$ is stable, and we have by Proposition \ref{st1} that the two-color stable multi-matching is a.s. unique and is attainable by the 2CIMC algorithm, we have that applying such algorithm, under $A$, to the coupled configuration a.s. leads to the graph which must be equal to $\tilde{B}$. This concludes the proof. $\square$\
 \
\section{Questions}
\begin{itemize}
 \item [(i)] Let $T$ be the total edge length of a typical point, that is, the sum of the length of all edges incident to it, and let $X$ be a probability distribution on the strictly positive integers. In \cite{Mia}, the following result was proved for the one-color case with points from a Poisson process $\mathcal{R}$ in $\mathbb{R}^{d}$.
\begin{theorem}(Theorem 1.1 \cite{Mia}).
There exists a translation invariant multi-matching scheme with $\mathbb{E}^*[T]<\infty$ if and only if $X$ has a finite moment of order $(d+1)/d$.
\end{theorem}
Consider the total edge length of a typical red point. Is there a similar result for two-color multi-matchings? 
\item [(ii)] Can the 2CIMC algorithm be analyzed when the marked point processes have infinite mean for the i.i.d. mark processes?
\item [(iii)] Give sharper conditions for percolation and non-percolation in the two-color stable multi-matching. In $d=1$, for the one-color case some advances have been made in \cite{Mia5}, when we assign a.s. 2 stubs per point. Could their methods and results be extended to the two-color case?
\end{itemize}

\textbf{Acknowledgements}: I thank my thesis advisor Mia Deijfen for many useful conversations and suggestions. I am also grateful to Alexander Holroyd who in my best knowledge first suggested the condition (\ref{condition}) of Theorem \ref{teo1}.

%\bibliography{abbr_long,pubext} 
\bibliographystyle{plain}
\bibliography{poissonb}

\begin{thebibliography}{10}

\bibitem{daley2}
Daryl~J. Daley and D.~Vere-Jones.
\newblock {\em An introduction to the theory of point processes. Vol II.
  General theory and structure}.
\newblock Springer-Verlag, New York, second edition, 2008.

\bibitem{Mia}
M.~Deijfen.
\newblock Stationary random graphs with prescribed iid degrees on a spatial
  {P}oisson process.
\newblock {\em Electronic Communications in Probability}, 14:81--89, 2009.

\bibitem{Mia2}
M.~Deijfen, O.~H\"aggstr\"om, and A.~E. Holroyd.
\newblock Percolation in invariant {P}oisson graphs with i.i.d degrees.
\newblock {\em Preprint, to appear in Arkiv for Matematik}.

\bibitem{Mia5}
M.~Deijfen, A.~E. Holroyd, and Y.~Peres.
\newblock Stable poisson graphs in one dimension.
\newblock {\em Electronic Journal of Probability}, 16:1238--1253, 2011.

\bibitem{Mia3}
M.~Deijfen and J.~Jonasson.
\newblock Stationary random graphs on $\mathbb{Z}$ with prescribed i.i.d.
  degrees and finite mean connections.
\newblock {\em Electronic Communications in Probability}, 11:336--346, 2006.

\bibitem{Mia4}
M.~Deijfen and R.~Meester.
\newblock Generating stationary random graphs on $\mathbb{Z}$ with prescribed
  i.i.d. degrees.
\newblock {\em Advances in Applied Probability}, 38(2):287--298, 2006.

\bibitem{Ferrari}
P.~A. Ferrari, C.~Landim, and H.~Thorisson.
\newblock Poisson trees, succession lines and coalescing random walks.
\newblock {\em Annales de l'Institut Henri Poincar\'e (B) Probabilit\'es et
  Statistiques}, 40(2):141--152, 2004.

\bibitem{gale}
D.~Gale and L.~Shapley.
\newblock College admissions and stability of marriage.
\newblock {\em The American Mathemathical Monthly}, 69, 1962.

\bibitem{hagg}
O.~H\"{a}ggstr\"om and Meester R.
\newblock Nearest neighbor and hard sphere models in continuum percolation.
\newblock {\em Random Structures \& Algorithms}, 9 (3):295--315, 1996.

\bibitem{Yuval1}
A.~E. Holroyd, R.~Pemantle, Y.~Peres, and O.~Schramm.
\newblock Poisson matching.
\newblock {\em Annales de l'Institut Henri Poincar\'e Probabilit\'es et
  Statistiques}, 45(1):266--287, 2009.

\bibitem{Yuval2}
A.~E. Holroyd and Y.~Peres.
\newblock Trees and matchings from point processes.
\newblock {\em Electronic Communications in Probability}, 8:17--27, 2003.

\bibitem{terry}
A.~E. Holroyd and T.~Soo.
\newblock Insertion and deletion tolerance of point processes.
\newblock {\em Preprint}, 2011.

\bibitem{jona}
J.~Jonasson.
\newblock Invariant random graphs with iid degrees in a general geography.
\newblock {\em Probability Theory and Related Fields}, 143:643--656, 2009.

\bibitem{kalle}
O.~Kallenberg.
\newblock {\em Foundations of modern probability. Probability and its
  Applications}.
\newblock Springer-Verlag, New York, second edition, 2002.

\bibitem{schon}
T.M. Liggett, R.H. Schonmann, and A.M. Stacey.
\newblock Domination by product measures.
\newblock {\em The Annals of Probability}, 25(1):71--95, 1997.

\end{thebibliography}
%\bibliography{poissonb}
\end{document}